\documentstyle[12pt]{article}
\begin{document}
\def\F(1){$ Fun(SU_{q}(2;\bf j))$}
\def\u2{\tilde u_2}
\def\z{\tilde z}
\def\t{\tilde }
\begin{center}
{\bf NONSEMISIMPLE  QUANTUM  GROUPS  AS  HOPF
ALGEBRAS  OF  THE  DUAL  FUNCTIONS}
\end{center}
\vspace{3mm}
\begin{center}
N.A. Gromov  \\
Department of Mathematics, Komi Research Center,
Ural Division, Russian Academy of Science,
Syktyvkar, 167000, Russia
\end{center}

{\bf Abstract.}
The  nonsemisimple quantum Cayley-Klein
groups \\ \F(1) are realized as Hopf algebra
of the noncommutative functions with the dual (or Study) variables.
The {\it dual} quantum algebras $ su_q(2;{\bf j}) $ are
constructed and their isomorphisms with the corresponding quantum
orthogonal algebras  $ so_q(3;{\bf j}) $ are established.
The possible couplings of the Cayley-Klein and Hopf
structures are considered.

\section{Introduction}
The simple (or semisimple) classical Lie groups (algebras)
may  be  trans\-form\-ed by contractions to the nonsemisimple
 ones
of the different structure. We named the set of such groups
(algebras) as Cayley-Klein (CK) groups (algebras).
In our approach
\cite{Mon--90} the contractions of the groups (algebras) are
described with the help of dual valued parameters, i.e.
the nonsemisimple CK groups are regarded as the groups over an
associative algebra $ D_{n}(\iota;\rm C) $ with the nilpotent
generators
 $ \iota_{k},
 \iota_{k}^2=0, k=1,...,n  $
 satisfying the commutative low of multiplications
$ {\iota_k}{\iota_m}={\iota_m}{\iota_k} \neq 0, k \neq m $.
The general element of the dual algebra $ D_{n}(\iota;\rm C) $
have the form
\begin{equation}
a=a_0 + \sum_{p=1}^{2^{n}-1} \sum_{k_1<...<k_p}
a_{k_1...k_p}\iota_{k_1}...\iota_{k_p}, \quad
a_0,a_{k_1...k_p}\in{\rm C}.
\end{equation}
For $ n=1 $ we have
$ D_{1}(\iota_{1};{\rm C})\ni{a=a_{0}+a_{1}\iota_{1}}, $
i.e. dual (or Study) numbers, when
$ a_{0},a_{1}\in{\rm R}. $
For $ n=2 $ the general element of
$ D_{2}(\iota_{1},\iota_{2};{\rm C}) $
is written as follows:
$ a=a_{0}+a_{1}\iota_{1}+a_{2}\iota_{2}+a_{12}\iota_{1}\iota_{2}. $


It is shown in \cite{Tr--94}, that the representations of
$ SU(2) $ , which are defined in the space of functions
on the qroup, are transformed to the representations of the
nonsemisimple euclidean qroup
 $ E(2), $ when the complex parameters of
 $ SU(2) $ are changed by the dual one's, i.e. when the
 complex variables of the functions are changed by the
 elements of the dual algebra
  $ D_{2}(\iota_{1},\iota_{2};{\rm C}). $

In present paper we apply this idea
to regard the nonsemisimple CK groups as the groups over the dual
algebra $ D_{n}(\iota;{\rm C}) $
to the case of quantum
groups. According with the general theory of quantum groups
\cite{Fad--89}, the quantum CK groups
\F(1) are regarded in Sec.2 as Hopf algebra of the
noncommutative functions of the dual variables.
The {\it dual}\footnote{To distinguish the different
meanings of the word "dual" we shall write
{\it dual} in the case of the dual aldebra.} to \F(1)
quantum algebras $  su_q(2;\bf j) $ are constructed in Sec.3
and their isomorphisms with the corresponding
quantum orthogonal CK algebras  $ so_q(3;{\bf j}) $
are found in Sec.4. The different dual parametrizations
of \F(1) lead to the different {\it dual}
algebras $ so_q(3;{\bf j}) $, which are distinguished by
the choise of the primitive element of the Hopf algebra. All
such possible combinations of the Hopf and CK structures are
considered in Sec.5.

\section{ Quantum group \F(1) }

The quantum unitary group $ Fun(SU_{\t q}(2)) $
is generated \cite{Fad--89} by the matrix with
non\-com\-mu\-ta\-ti\-ve elements
      $$
T= \left (  \begin{array}{cc}
      {\t a}  & {\t b}  \\
      -{\t q^{-1}}\bar {\t b}  &  \bar {\t a}   \end{array}  \right ) =
      \left ( \begin{array}{cc}
      {\t a}_1+i{\t a}_2    &    {\t b}_1+i{\t b}_2 \\
      -e^{-\t z}({\t b}_1-i{\t b}_2)  & {\t a}_1-i{\t a}_2
 \end{array}  \right ),
$$
\begin{equation}
{\det}_qT={\t a}_1^2+{\t a}_2^2+{\t b}_1^2+ {\t b}_2^2=1,
\label{1.1}
\end{equation}
where
the deformation
 parameter $ \t q=e^{\t z}, \;
 {\t z}\in{\rm C} $ and the bar designate the complex
conjugation. Commutation relations for generator are defined by the
following R-matrix
\begin{equation}
R_{\t z}=\left (  \begin{array}{cccc}
      {\t q} & 0 & 0 & 0 \\
      0 & 1 & 0 & 0 \\
      0 & \t {\lambda} & 1 & 0 \\
      0 & 0 & 0 & {\t q}  \end{array} \right )=
      \left ( \begin{array}{cccc}
      e^{\t z} & 0 & 0 & 0 \\
      0 & 1 & 0 & 0 \\
      0 & 2\sinh{\t z} & 1 & 0 \\
      0 & 0 & 0 & e^{\t z} \end{array} \right ),
\label{1.2}
\end{equation}
where $ \t {\lambda}=\t q-\t q^{-1}=2\sinh{\t z}. $

 {\bf Proposition 1.} {\it The quantum group} \F(1),
 $ {\bf j}=(j_1,j_2), \; j_k=1,\iota_k, \; k=1,2 $
      {\it is given by }
      $$
T({\bf  j})=
 \left ( \begin{array}{cc}
      a    &    b \\
      -e^{-j_1j_2z}{\bar b} & \bar a \end{array} \right ) =
\left ( \begin{array}{cc}
 a_1+ij_1j_2a_2               & j_1b_1+ij_2b_2 \\
-e^{-j_1j_2z}(j_1b_1-ij_2b_2) &  a_1-ij_1j_2a_2
 \end{array} \right ),
      $$
\begin{equation}
{\det}_qT({\bf j})=a_1^2+j_1^2j_2^2a_2^2+j_1^2b_1^2+j_2^2b_2^2=1,
\label{1.3}
\end{equation}
$$
{R_z({\bf j})}{T_1({\bf j})}{T_2({\bf j})}=
{T_2({\bf j})}{T_1({\bf j})}{R_z({\bf j})},
$$
\begin{equation}
T_1({\bf j})=T({\bf j})\otimes I, \quad
T_2({\bf j})=I\otimes T({\bf j}),
\label{1.6}
\end{equation}
\begin{equation}
\triangle T({\bf j})=T({\bf j})\dot {\otimes}T({\bf j}),\quad
\epsilon(T({\bf j}))=I,
\label{1.7}
\end{equation}
\begin{equation}
S(T({\bf j}))= \left ( \begin{array}{cc}
 a_1-ij_1j_2a_2   &  -e^{j_1j_2z}(j_1b_1+ij_2b_2) \\
 j_1b_1-ij_2b_2   &  a_1+ij_1j_2a_2
 \end{array}   \right ),
\label{1.8}
\end{equation}
{\it where}
\begin{equation}
R_{z}({\bf j}) =\left (  \begin{array}{cccc}
      e^{j_1j_2z} & 0 & 0 & 0 \\
      0 & 1 & 0 & 0 \\
      0 & 2\sinh{j_1j_2z}& 1 & 0 \\
      0 & 0 & 0 & e^{j_1j_2z} \end{array} \right ).
\label{1.5}
\end{equation}

In explicit form the commutation relation are as follows
      $$
[b_1,b_2]=0, \quad
 [a_1,a_2]=
 -i(j_1^2b_1^2+j_2^2b_2^2)e^{j_1j_2z}j_1^{-1}j_2^{-1}\sinh{j_1
 j_2z},
      $$
\begin{eqnarray}
 a_1b_1                 &                 =                 &
 b_1a_1\cosh{j_1j_2z}+ij_1j_2b_1a_2\sinh{j_1j_2z}, \nonumber \\
 a_1b_2                 &                 =                 &
 b_2a_1\cosh{j_1j_2z}+ij_1j_2b_2a_2\sinh{j_1j_2z}, \nonumber \\
 a_2b_1                 &                 =                 &
 b_1a_2\cosh{j_1j_2z}-ib_1a_1j_1^{-1}j_2^{-1}\sinh{j_1j_2z},
 \nonumber \\
 a_2b_2 & = &
 b_2a_2\cosh{j_1j_2z}-ib_2a_1j_1^{-1}j_2^{-1}\sinh{j_1j_2z}.
 \label{1.9}
\end{eqnarray}
 In particular for the most contracted case
 $ j_1=\iota_1, j_2=\iota_2 $
 we have from Eq.(\ref{1.3})
      $ a_1=1 $ and Eqs.(\ref{1.9}) become as follows:
 \begin{equation}
  [b_1,b_2]=0, \; [b_1,a_2]=izb_1, \; [b_2,a_2]=izb_2.
  \label{1.10}
 \end{equation}

  The matrix (\ref{1.3}) and R-matrix (\ref{1.5}) are
obtained from (\ref{1.1}) and (\ref{1.2})
by the following transformations of generators and
deformation parameter:
\begin{equation}
\t z=j_1j_2z, \; \t a_1=a_1, \; \t a_2=j_1j_2a_2, \;
 \t b_1=j_1b_1, \; \t b_2=j_2b_2.
\label{1.4}
\end{equation}

\section{ $ su_q(2;{\bf j}) $ as the dual to
 \F(1) }

 The generators of the {\it dual} to $ Fun(SU_{\tilde q}(2)) $
 quantum algebra
 $ su_{\tilde q}(2) $ are written \cite{Fad--89} in compact matrix form as
\begin{equation}
L^{(+)}= \left (  \begin{array}{cc}
  \tilde t   &   \tilde u_1+i\u2 \\
    0   &  \tilde t^{-1}   \end{array}   \right ), \qquad
L^{(-)}=
    \left (  \begin{array}{cc}
    \tilde t^{-1}   &    0    \\
      -e^{\z}(\tilde u_1-i\u2) & \tilde t    \end{array} \right  ).
\label{2.1}
\end{equation}

Following \cite{Fad--89}, we define the quantum
algebra $ su_q(2;\bf j),  $  {\it dual} to \F(1) as
\begin{equation}
<L^{(\pm)}({\bf j}),T({\bf j})>=R^{(\pm)}({\bf j}),
\label{2.2}
\end{equation}
where $ L^{(\pm)}({\bf j}) $ are given by
      $$
L^{(+)}(\bf j)=\left (    \begin{array}{cc}
           t    &  u  \\
           0    &  t^{-1}  \end{array}  \right  )=
   \left  (  \begin{array}{cc}
      t  &  j_1^{-1}u_1+ij_2^{-1}u_2 \\
      0  &   t^{-1}   \end{array}  \right  ),
      $$
\begin{equation}
L^{(-)}({\bf j})= \left (  \begin{array}{cc}
  t^{-1}   & 0  \\
  -e^{j_1j_2z}\bar u &  t  \end{array}  \right  )=
  \left  (    \begin{array}{cc}
  t^{-1}  &   0  \\
  -e^{j_1j_2z}(j_1^{-1}u_1-ij_2^{-1}u_2)  &    t  \end{array} \right ),
 \label{2.3}
  \end{equation}
and
act on the first order polinomial of the generators of \F(1).
 $ R^{(\pm)}({\bf j}) $ are expressed by R-matrix
(\ref{1.5}) as
\begin{eqnarray}
R({\bf j})=e^{-j_1j_2z}R_z({\bf j}),& \det R({\bf j})=1,& \nonumber \\
R^{(-)}({\bf j})=R^{-1}({\bf j}),&  R^{(+)}({\bf j})=PR({\bf j})P,&
\nonumber \\
P(a\otimes b)=b\otimes a.& &
\label{2.4}
\end{eqnarray}
The matrices (\ref{2.3}) may be obtained from (\ref{2.1})
by the following (contraction) transformations of the generators and
the deformation parameter
\begin{equation}
\t z=j_1j_2z, \quad \tilde {u_1}=u_1/j_1, \quad \u2=u_2/{j_2},\quad
\t H=H/{j_1j_2},
\label{2.5}
\end{equation}
where $ \t t=\exp({\t z}{\t H}/2), t=\exp(zH/2). $
In explicit form the  actions (\ref{2.2}) are given by
      $$
 t(b)=t(\bar b)=u(b)=\bar u(\bar b)=u(a)=u(\bar a)=0,
      $$
\begin{equation}
t(a)=x, \quad t(\bar a)=x^{-1},\quad
u(\bar b)=-x\lambda, \quad \bar u(b)=x^{-1}\lambda,
\label{2.6}
\end{equation}
where  $ x=e^{j_1j_2z/2}, \; \lambda=2\sinh{j_1j_2z}. $
Equations (\ref{2.6}) become
      $$
 u_k(a_1)=u_k(a_2)=0, \; t(a_1)= \cosh j_1j_2z/2, \;
 t(a_2)=-ij_1^{-1}j_2^{-1}\sinh{j_{1}j_{2}z/2},
      $$
\begin{eqnarray}
u_k(b_k) & = & -\sinh{j_1j_2z} \cdot \sinh{j_1j_2z/2},
 \; k=1,2, \nonumber \\
u_1(b_2) & = & -ij_1^2{{\sinh{j_1j_2z}}\over {j_1j_2}}
 \cdot{\cosh{j_1j_2z/2}}, \nonumber \\
u_2(b_1)   &   =   &  ij_2^2{{\sinh{j_1j_2z}}\over  {j_1j_2}}
 \cdot{\cosh{j_1j_2z/2}} .
\label{2.7}
\end{eqnarray}

   {\it  Remark  1.}  Eqs.(\ref{2.7})  for  $  j_1=j_2=1 $
 describe the quantum algebra $ su_q(2). $

{\it  Remark 2.} At first sight nondiagonal elements of $ L^{(\pm)}({\bf j}) $
(\ref{2.3}) are not defined for dual values of $ j_k $,
since the division of a real or complex number by the dual units
 $ \iota_k $ are not defined. But the matrices $ L^{(\pm)}({\bf j}) $
are a {\it linear} functionals on the dual variables
 $ j_kb_k  $, therefore their  {\it actions}
on the generators of \F(1) gives the well defined
expressions  (\ref{2.7}).

 {\bf Proposition 2.} {\it The quantum algebra} $ su_q(2;{\bf j}) $
 {\it is given by}
\begin{equation}
R^{(+)}({\bf j})L_1^{(+)}({\bf j})L_2^{(-)}({\bf j})=
L_2^{(-)}({\bf j})L_1^{(+)}({\bf j})R^{(+)}({\bf j}),
\label{2.8}
\end{equation}
\begin{eqnarray}
u_{1}t & = & tu_1\cosh{j_{1}j_{2}z}+ij_{1}^{2}tu_{2} {{\sinh{j_{1}j_{2}z}}
\over {j_{1}j_{2}}},
 \nonumber  \\
u_{2}t & = & tu_2\cosh{j_{1}j_{2}z}-ij_{2}^{2}tu_{1} {{\sinh{j_{1}j_{2}z}}
\over {j_{1}j_{2}}}, \nonumber
\end{eqnarray}
\begin{equation}
[j_{1}^{-1}u_{1}, j_{2}^{-1}u_{2}]  =
 -2ie^{-j_1j_2z}\cdot {\sinh{j_1j_2z}}\cdot {\sinh{zH}}.
\label{2.9}
\end{equation}
$$
\triangle L^{(\pm)}({\bf j})=L^{(\pm)}({\bf j}){\dot \otimes}L^{(\pm)}({\bf j})
,
 \quad
\epsilon(L^{(\pm)}({\bf j}))=I,
$$
      $$
\triangle  t = t \otimes t, \quad
\triangle  u_k  =  t  \otimes u_k + u_k \otimes t^{-1}, \quad
 k=1,2,
      $$
\begin{equation}
\epsilon(t)  =  1, \quad \epsilon(u_1)=\epsilon(u_2)=0,
\label{2.11}
\end{equation}
$$
S(L^{(+)}({\bf j}))= \left (  \begin{array}{ll}
  t^{-1} &  -e^{j_1j_2z}(j_1^{-1}u_1+ij_2^{-1}u_2) \\
  0    &   t   \end{array}  \right  ),
  $$
\begin{equation}
S(L^{(-)}({\bf j}))= \left (  \begin{array}{ll}
  t      &   0  \\
(j_1^{-1}u_1-ij_2^{-1}u_2)      &   t^{-1}   \end{array}  \right  ).
\label{2.12}
\end{equation}

 \section{ Isomorphism of $ su_q(2;{\bf j}) $ and $ so_q(3;{\bf j})$ }

We have described the quantum algebra $ su_q(2;{\bf j}) $
as the {\it dual} $ \; $ to \F(1). In this section we will show their
isomorphism with the orthogonal CK algebra
 $ so_q(3;{\bf j}). $ The quantum analogue of the  universal enveloping
 algebra of $ so(3;{\bf j}) $ with
 the rotation generator $ X_{02} $ is the primitive
element has been given in \cite{Val--93}:

 {\bf Proposition 3.} {\it
 The Hopf algebra structure of}
 $ so_q(3;{\bf j};X_{02}) $ {\it is given by}
      $$
\triangle X_{02}  =  I\otimes X_{02}+X_{02}\otimes I,
      $$
      $$
\triangle X  =  e^{-{\hat z}X_{02}/2}\otimes X+X\otimes e^{{\hat z}X_{02}/2},
\quad  X=X_{01},X_{12},
      $$
      $$
\epsilon(X_{01})  =  \epsilon(X_{02})=\epsilon(X_{12})=0, \quad
\t S(X_{02})  =  -X_{02},
      $$
      $$
\t S(X_{01})  =  -X_{01}\cos{j_1j_2{\hat z}/2}+
j_1^2X_{12}{{\sin{j_1j_2{\hat z}/2}\over {j_1j_2}}},
      $$
      $$
\t S(X_{12})  =  -X_{12}\cos{j_1j_2{\hat z}/2}-
j_2^2X_{01}{{\sin{j_1j_2{\hat z}/2}\over {j_1j_2}}},
      $$
 \begin{equation}
 [X_{01},X_{02}]=j_1^2X_{12},\quad   [X_{02},X_{12}]=j_2^2X_{01},\quad
 [X_{12},X_{01}]={\sinh({\hat z}X_{02})\over {\hat z}}.
\label{3.2}
 \end{equation}

 {\bf Proposition 4.} {\it The Hopf algebra} $ su_q(2;{\bf j}) $
      (\ref{2.8})--(\ref{2.12})  {\it  is  isomorphic to the
Hopf algebra} $ so_q(3;{\bf j};X_{02}). $

      {\it Proof.}
It is easily verify, that the equations of the previous
section are passed in the corresponding expression  (\ref{3.2})
by the following substitution of the deformation parameter $ z $
and the generators $ t, j_k^{-1}u_k $ for the deformation
parameter $ {\hat z} $ and the generators $ X $:
      $$
z  =  i{\hat  z}/2,  \quad  t=e^{zH/2},  \quad   H=-2iX_{02},
      $$
        $$
 j_1^{-1}u_1  = 2ij_1De^{-ij_1j_2{\hat z}/4}X_{12}, \quad
j_2^{-1}u_2  =  2ij_2De^{-ij_1j_2{\hat z}/4}X_{01},
        $$
 \begin{equation}
D  =  i{\left ( {{\hat z}\over {2j_1j_2}}\sin(j_1j_2{\hat z}/2)
 \right ) }^{1/2}.
\label{3.3}
\end{equation}
The different definitions of the antipode $ \t S $ in \cite{Val--93} as
\begin{equation}
\t S(X)=-e^{{\hat z}X_{02}/2}Xe^{-{\hat z}X_{02}/2}=
-e^{ zH/2}Xe^{-zH/2},
\label{3.4}
\end{equation}
and the antipode (\ref{2.12}) of $ su_q(2;{\bf j}) $ as
\begin{equation}
S(u_k)=-e^{-zH/2}u_k{e^{zH/2}},
\label{3.5}
\end{equation}
must be taken into account. This produce
a slight difference of signs when Eqs. (\ref{2.12})
are trans\-for\-med by
 (\ref{3.3}) as compared with Eqs. (\ref{3.2}).

\section{On the possible couplings of Hopf and \-
 Cayley--Klein structures}

The quantum CK algebra $ so_q(n+1;{\bf j}) $ can be
got starting from the quantum orthogonal algebra
 $ so_q(n+1) $ by the following (contraction) transformations of the
 generators and deformation parameter
\begin{equation}
X_{\mu \nu}  =  J_{\mu \nu}{\t X_{\mu \nu}},  \quad
J_{\mu \nu}=\prod^{\nu}_{l=\mu+1}{j_l}, \quad
\mu<\nu, \quad
\t z  =  Jz ,
\label{4.1}
\end{equation}
where the explicit expression for
multiplier $ J $ depends on the choose of the primitive
elements of the Hopf algebra. Note that for nonquantum case
the transformations of generators (\ref{4.1}) give the CK
algebra
$  so(n+1;{\bf j}) $ \cite{Mon--90}, \cite{TMF--81},  \cite{JMP--92}.

Under introduction of a Hopf algebra structure in the
universal enveloping algebra some commuting generators
of $ so_q(n+1) $ ( the basic ones in the Cartan subalgebra)
have to be selected as the primitive elements of the
Hopf algebra, i.e. they are
distinguished among other generators.
Let us observe, that any permutation
 $ \sigma \in S(n+1) $  indices of the rotation generators
$ \t X_{\mu \nu} $ of $ so_q(n+1) $
(i.e. the transformations from the Weyl group)
leads to the isomorphic Hopf algebra
$ so_q^{'}(n+1), $ but with some other set of primitive
generators. This isomorphism of the Hopf algebras may be
destroyed by contractions. Indeed, under transformations
(\ref{4.1}) the primitive generators of
$ so_q(n+1) $ and $ so_q^{'}(n+1) $ are multiplied by
the different
 $ J_{\mu \nu} $ and for the specific dual values of
 $ j_k, $ this leads to the nonisomorphic quantum algebras
 $ so_q(n+1;{\bf j}) $ and
$ so_q^{'}(n+1;{\bf j}) $. In other words, the Hopf and CK
structures may be combined in a different manner
 for quantum CK algebras.

Let us illustrate the different couplings of Hopf and
CK structures for the case of
 $ so_q(3;{\bf j}) $.
The Hopf algebra $ so_q(3; \t X_{02}) $ with
$ \t X_{02} $ as the primitive element is described by
 Eqs. (\ref{3.2})
for $ j_1=j_2=1 $. Transformations (\ref{4.1})
are in this case as follows:
\begin{equation}
X_{01}=j_1{\t X_{01}},\quad X_{02}=j_1j_2{\t X_{02}},\quad
X_{12}=j_2{\t X_{12}},\quad \t {\hat z}=j_1j_2{\hat z}
\label{4.2}
\end{equation}
and give in result Eqs. (\ref{3.2}).
The transformation law of the deformation parameter
is defined by those of the primitive generator
 $ \t X_{02} $, i.e. $ J=j_1j_2.$
The substitution $ \sigma \in S(3),\; \sigma(0)=1, \sigma(1)=0, \sigma(2)=2  $
indices of the generators
 $ \t X_{\mu \nu} $ leads to the Hopf algebra $ so_q(3; \t X_{12}) $
 with $ \t X_{12} $ as the primitive element. The introduction
 of the CK structure by the transformations of generators
as in Eqs.(\ref{4.2}) and the deformation parameter like
 $ \t X_{12}, $ i.e.
 $  \t  {\hat  z}=j_2{\hat  z}  $  leads  to  the   following
 proposition.

 {\bf Proposition 5.} {\it
 The quantum algebra} $ so_q(3;{\bf j}; X_{12}) $
     {\it is given by }
$$
\triangle X_{12}  =  I\otimes X_{12}+X_{12}\otimes I,
$$
$$
\triangle X  =
 e^{-X_{12}{\hat z}/2} \otimes X+X \otimes e^{X_{12}{\hat z}/2},
\quad X=X_{01}, X_{02},
$$
$$
\epsilon(X_{01})  =  \epsilon(X_{02})=\epsilon(X_{12})=0, \quad
{\t S(X_{12})} =-X_{12},
$$
$$
{\t S(X_{01})}  =
 -X_{01}\cos{j_2{\hat z}/2}-X_{02}j_2^{-1}\sin{j_2{\hat z}/2},
$$
$$
{\t S(X_{02})}  =  -X_{02}\cos{j_2{\hat z}/2}+X_{01}j_2\sin{j_2{\hat z}/2},
$$
\begin{equation}
{[X_{01},X_{02}]}=j_1^2{\hat z}^{-1}\sinh({\hat z}X_{12}),\;
{[X_{02},X_{12}]}=j_2^2X_{01},\; {[X_{12},X_{01}]}=X_{02}
\label{4.3}
\end{equation}

      {\bf Proposition 6.} {\it
The  quantum qroup} \F(1) {\it dual to} (\ref{4.3})
{\it is given by relations} (\ref{1.6}),(\ref{1.7}) {\it with }
$$
T({\bf j})=
      \left ( \begin{array}{cc}
       a_1+ij_2a_2   &    j_1(b_1+ij_2 b_2) \\
      -e^{-j_2 z}j_1(b_1-ij_2 b_2)  & a_1-ij_2a_2    \end{array}  \right ),
$$
\begin{equation}
{\det}_q{T({\bf j})}=a_1^2+j_2^2a_2^2+j_1^2b_1^2+j_1^2j_2^2b_2^2=1,
\label{4.4}
\end{equation}

 {\it and}

\begin{equation}
R_{z}({\bf j}) =\left (  \begin{array}{cccc}
      e^{j_2z} & 0 & 0 & 0 \\
      0 & 1 & 0 & 0 \\
      0 & 2\sinh{j_2z}& 1 & 0 \\
      0 & 0 & 0 & e^{j_2z} \end{array} \right ).
\label{4.4'}
\end{equation}

      The matrices (\ref{4.4}),(\ref{4.4'})
 are obtained from (\ref{1.1}),(\ref{1.2}) by transformations
\begin{equation}
\t z=j_2z, \; \t a_1=a_1, \; \t a_2=j_2a_2, \;
 \t {b_1}=j_1b_1, \; \t {b_2}=j_1j_2b_2.
\label{4.5}
\end{equation}

     The {\it dual} to (\ref{4.4}) quantum algebra is defined
 by Eqs.(\ref{2.2}), where
the matrices $ L^{(\pm)}({\bf j}) $ are as follows
      $$
L^{(+)}({\bf j})= \left (  \begin{array}{cc}
    t   &   j_1^{-1}(u_1+ij_2^{-1}u_2) \\
    0   &   t^{-1}   \end{array}   \right ),
      $$
\begin{equation}
L^{(-)}({\bf j})=
    \left (  \begin{array}{cc}
      t^{-1}   &    0    \\
      -e^{j_2z}j_1^{-1}(u_1-ij_2^{-1}u_2) &  t    \end{array} \right  )
\label{4.7}
\end{equation}
and are obtained from (\ref{2.1}) by (contraction) transformations
\begin{equation}
\t z=j_2z, \quad \t H=H/j_2, \quad \t u_1=u_1/j_1,
\quad \t u_2=u_2/j_1j_2,
\label{4.8}
\end{equation}
where $ t=\exp(zH/2) $. Nonzero actions of
$ t, u_k $ on $ a, b_k $ are found from the relations
 (\ref{2.2}) in the form
      $$
 t(a_1)= \cosh{j_{2}z/2}, \quad
 t(a_2)= -ij_2^{-1}\sinh{j_{2}z/2},
      $$
      $$
 u_1(b_1) = u_2(b_2) = -\sinh{j_2z}\cdot \sinh{j_2z/2},
      $$
      $$
u_1(b_2)  =  -ij_2^{-1}\sinh{j_2z}\cdot \cosh{j_2z/2},
      $$
\begin{equation}
u_2(b_1)  =  ij_2\sinh{j_2z}\cdot \cosh{j_2z/2},
\label{4.9}
\end{equation}
and commutators follow from Eqs.(\ref{2.8})
      $$
[H,u_1] = -2iu_2, \quad [H,u_2]=2ij_2^2u_1 ,
      $$
\begin{equation}
[u_2,u_1] = 2ij_1^2j_2e^{-j_2z}\sinh{j_2z}\cdot \sinh{zH}.
\label{4.10}
\end{equation}

      {\bf Proposition 7.} {\it
      Hopf algebra} (\ref{4.7})--(\ref{4.10}) {\it is isomorphic to }
$ so_q(3;{\bf j};X_{12}) $ (\ref{4.3}).

      {\it Proof.}
 The connections of $ H, u_k $ with the rotation generators $ X_{\mu\nu} $
of $ so_q(3;{\bf j};X_{12}) $ are given by
$$
z=i{\hat z}/2, \quad H=-2iX_{12},\quad u_1=F\cdot X_{02},
\quad j_2^{-1}u_2=-j_2F\cdot X_{01},
$$
\begin{equation}
F=e^{-ij_2{\hat z}/4}{(2{\hat z}j_2^{-1}\sin(j_2{\hat z}/2))}^{1/2}.
\label{4.11}
\end{equation}

As it was mention early the contractions of quantum group
 and algebras are correspond
to the dual values of the parameters $ j_k. $
In particular, the quantum euclidean algebra  $ so_q(3;\iota_1,1;X_{12}) $
is described by (\ref{4.3}) for  $ j_1=\iota_1, j_2=1. $
The deformation parameter is left untoched, since
 $ j_2=1 $.
Their {\it dual} euclidean quantum group is
realized according with
(\ref{4.4}) as Hopf algebra of the noncommutative
functions with dual variables
(cf. \cite{Vak--89}--\cite{Mas--94}).

The third possible coupling of Hopf and CK structures
is connected with the choose of $ \t X_{01} $ as the
primitive element and leads to the quantum algebra $ so_q(3; \t X_{01}) $.

      {\bf Proposition 8.} {\it
The  quantum CK algebra} $ so_q(3;{\bf j};  X_{01}) $
 {\it is given by }
      $$
\triangle X_{01}  =  I\otimes X_{01}+X_{01}\otimes I,
      $$
      $$
\triangle X  =
 e^{-X_{01}{\hat z}/2} \otimes X + X \otimes e^{X_{01}{\hat z}/2},
\quad X=X_{02}, X_{12},
      $$
      $$
\epsilon(X_{01})  =  \epsilon(X_{02})=\epsilon(X_{12})=0, \quad
\t S(X_{01})=-X_{01},
      $$
      $$
{\t S(X_{02})} = -X_{02}\cos{j_1{\hat z}/2}-X_{12}j_1\sin{j_1{\hat z}/2},
      $$
      $$
{\t S(X_{12})} = -X_{12}\cos{j_1{\hat z}/2}+X_{02}j_1^{-1}\sin{j_1{\hat z}/2},
      $$
\begin{equation}
{[X_{01},X_{02}]}  =  j_{1}^{2}X_{12}, \;
{[X_{12},X_{01}]}=X_{02}, \;
 {[X_{02},X_{12}]}=j_2^2{{\sinh{{\hat z}X_{01}}}\over {\hat z}}.
\label{4.12}
\end{equation}

Last equations are obtained in a three steps:
i) put $ j_1=j_2=1 $ in Eqs. (\ref{3.2}),
ii) apply permutation  $ \sigma \in S(3),  \;
     \sigma(0)=0,\,\sigma(1)=2,\,\sigma(2)=1 $,
 iii) introduce CK structure by transformation
(\ref{4.2}) of generators and deformation parameter like
$ \t X_{01} $, i.e.  $ {\t {\hat z}}=j_1{\hat z}. $
The {\it dual} to $ so_q(3;{\bf j};X_{01}) $
quantum group \F(1) is obtained from (\ref{1.1})
by the transformations
\begin{equation}
\t z=j_1z, \; \t a_1=a_1, \; \t a_2=j_1a_2, \;
 \t b_1=j_1j_2b_1, \; \t b_2=j_2b_2,
\label{4.13}
\end{equation}
and the matrices $ L^{(\pm)}({\bf j}) $ are obtained from (\ref{2.1})
by transformations
\begin{equation}
\t z=j_1z, \quad \t H=H/j_1, \quad \t u_1=u_1/{j_1j_2}, \quad
\t u_2=u_2/j_2.
\label{4.14}
\end{equation}

Eqs. (\ref{4.12}) for $ j_1=\iota_1, j_2=1 $
describe the quantum euclidean algebra, which has been
obtained in \cite{Cel--90} by contraction of
 $ su_q(2) $. The deformation parameter is transformed
 in this case.

Eqs. (\ref{3.2}) and (\ref{4.3}) for
$ j_1=\iota_1, j_2=\iota_2 $
give two quantum galilean algebras, which has been obtained in
\cite{Val--93}, \cite{Cel--90}.
(Galilean algebra (\ref{4.12}) is isomorphic to the algebra (\ref{4.3})
for $ j_1=\iota_1, j_2=\iota_2 $).

\section{ Conclusion}

Contractions are the method of receiving a new Lie groups (algebras)
from the unitial ones, in particular, the nonsemisimple qroups
(algebras) from the simple ones. In the traditional approach
\cite{Wig--53} this is achieved by introduction of a real zero
tending parameter $ \epsilon \rightarrow 0 $. In our approach
\cite{Mon--90},  \cite{TMF--81},  \cite{JMP--92}
contractions are described by the dual valued parameters
 $ j_k $.
In the case of standart matrix theory of quantum
groups
these contractions supplemented with the
appropriate transformations of the deformation parameter
lead to realization of nonsemisimple quantum qroups as Hopf
 algebras of noncommutative functions with dual variables.
 Although both descriptions of contractions are equivalent in many
 respects, it seems that the language of dual (or Study) numbers
          is mathematically more correct.
 For example, the quantum euclidean qroup is described by the matrix
 (\ref{4.4}) for
$ j_1={\iota}_1, j_2=1  $ with dual nondiagonal elements, while
the limit
 $ l \rightarrow 0 $ in the traditional approach gives the matrix
 with zero nondiagonal elements
 (cf. (16) in \cite{Zum--92}).
 The dual units are nilpotent like the Grassmanian ones, which
 are used for description of supersymmetry. The only difference
 is the commutative or anticommutative law of multiplications.
 The dual algebra $ D_{n}(\iota;{\rm C}) $ is the subalgebra of the
 even part of the Grassmanian algebra with $ 2n $ nilpotent
 generators.

The constructive algorithm for the description of
a different couplings of Hopf and CK structures
is given.
The different combinations of Hopf and CK structures for the
quantum algebras $ so_q(3;{\bf j}) $ are appear on the level
of the quantum groups \F(1) as the different dual functions
for the  elements of matrix $ T $
(cf.(\ref{1.3}),(\ref{4.4})).
The  transformation  law  of the deformation parameter is the
 same as the one of the dual part of the diagonal elements of
$ T $ (cf.(\ref{1.4}),(\ref{4.5}),(\ref{4.14})).

\section{ Acknowledgements }

I would like to thank A.Ballesteros, F.J.Herranz, M.A.del Olmo,
 M.San\-tan\-der and I.V.Kostyakov for usefull discussions.
This work has been partially supported by the Russian
Foundation of Fundamental Research (project 93--011--16180).

\end{document}